\documentclass[12pt]{article}
\usepackage{geometry}
\usepackage{amssymb}
\usepackage{color}

\usepackage{amsmath}

\newif\ifleft
\lefttrue % \*true or \*false
\newif\ifpre
\pretrue % \*true or \*false
\newif\ifhide
\hidefalse

\def\1{\mathbf{1}}

\newtheorem{theorem} {Theorem}

\newtheorem{remark} {Remark}

\newif \ifshowup
\showupfalse %true or false

\title{{\normalsize\tt\hfill\jobname.tex}\\
      \bf On pathwise uniqueness for  
      multidimensional McKean--Vlasov equations 
      \\~\\
A. Yu. Veretennikov\footnote{University of Leeds, UK
$\,$ \& $\,$ National Research University Higher School of Economics, and Institute for Information Transmission Problems, Moscow, Russian Federation; email: a.veretennikov@ leeds.ac.uk}
}

\begin{document}
\maketitle

%{\color{magenta}(to make 13 pages! e.g., 1.5 intro + 2.5 main results + 2.5 refs = 6.5; remains up to 6.5 pages for the sketches of the proofs)}

\begin{abstract}
Pathwise uniqueness  for multi-dimensional stochastic McKean--Vlasov equation is established under moderate  regularity conditions on the drift and diffusion coefficients. Both drift and diffusion depend on the marginal measure of the solution. For pathwise uniqueness, the drift is assumed to be Dini-continuous in the state variable, while the diffusion must be Lipschitz, continuous in time and uniformly nondegenerate. The setting is classical McKean--Vlasov, that is, coefficients of the equation are represented as integrals over the marginal distributions of the process. 
\end{abstract}

\section{Introduction}
%\subsection{Setting, backgrounds and motivation (sokratit'!)}
Consider solutions of a stochastic McKean--Vlasov equation in $\mathbb R^d$
\begin{equation}\label{e1}
dX_t = B[t,X_t, \mu_t]dt + \Sigma[t,X_t, \mu_t]dW_t, \qquad
X_0=x_0,
\end{equation}
with a possibly random $x_0$ independent on $W$ (and possessing a second moment), under the convention
\begin{equation}\label{e200}
B[t,x,\mu]=\int
b(t,x,y)\mu(dy), \;\; \Sigma[t,x,\mu]=\int
\sigma(t,x,y)\mu(dy).
\end{equation}
Here $W$ is a standard 
$d_1$-dimensional
Wiener process with $d_1\ge d$, $b$ and $\sigma$ are vector and matrix Borel
functions of corresponding dimensions $d$ and $d\times d_1$,  respectively, $\mu_t$ is the distribution of the process $X$ at
$t$. The initial data $x_0$ may be random, but independent of
$W$; a delta-measure is allowed. The systematic study of such equations
was started by McKean \cite{McKean}. The reference \cite{Sz} provides an introduction to the whole area. 
McKean--Vlasov's equations are interesting and important in quite a few areas such as 
multi-agent SDE systems, 
%(see \cite{Ben, Bossy_Talay}), 
filtering,  
%(see \cite{Crisan_Xiong}),  
et al. 

Results on strong solutions for the equation (\ref{e1}) can be found, in particular, in \cite{Sz}, \cite{Ver06}, with some last achievements in \cite{Bahlali19} and \cite{Paul-Eric}. About most recent results for diffusion equations with jumps see also \cite{SimaMehri2}. We suggest another method in comparison to \cite{Paul-Eric} with the goal to prove some result by using a combination of the approaches from \cite{Bahlali19}, \cite{Zvonkin} and \cite{YamadaWatanabe}. The point is to use  relaxed assumptions on the dependence of the drift  coefficient with respect to the state variable in comparison to \cite{Paul-Eric},   %(Dini conditions instead of Hoelder), 
to avoid use differentiation of measures, and as a result to have a significantly shorter presentation of the result.  
Yet, the assumptions on the ``third variable'' in our paper are stronger than in \cite{Paul-Eric}, so that, strictly speaking, the overall sets of conditions in this and in earlier papers are not directly comparable. 

This short paper consists of the introduction in the section 1,  the main result in the section 2, and its proof in the section 3.

\section{Main result}\label{sec:we}

Note that for any Borel function $f(z, y)$ and any probability measure
$\mu(dy)$ such that $f(z, \cdot)$ is integrable with respect to this measure, the
function $\displaystyle f[z, \mu]:= \int f(z,y)\,\mu(dy)$ is   Borel measurable in $z$ (e.g., cf. \cite[Theorem 2.6.8]{Shi}).
So, in particular, 
if for each $(t,x)$ the Borel measurable coefficients $b(t,x,y)$ and $\sigma(t,x,y)$ are bounded -- which we have assumed-- then the functions $\tilde b(t,x):=
B[t, x, \mu_t]$ and $\tilde \sigma(t,x):=
\Sigma[t, x, \mu_t]$ are Borel measurable in $(t, x)$. Due to this fact, the equation (\ref{e1}) is well-posed. Denote by ${\cal P}_2$ the set of all probability measures in ${\mathbb R}^d$ with a finite second moment.

%Assume {\color{red}(mozhet, ne nado dlia edinstvenoosti? a rabotat' s uslovnym m.o.?)}
%{\color{red} too much???}
%\begin{equation}\label{m4}
%{\mathbb E}|x_0|^2<\infty.
%\end{equation} 
%and  weak existence for the equation (\ref{e1}) (cf. ???).
Denote
 $$
 \rho_B(r):= \sup_{t\ge 0}\, \sup_{\mu \in {\cal P}_2}\sup_{|x-x'|\le r}|B[t,x,\mu] - B[t,x',\mu]|, \quad r\ge 0,
$$ 
which is a ``uniform with respect to other variables'' modulus of continuity of the function $B$ in the $x$ variable. The notation ${\cal L}(X_t)$ stands for the marginal distribution of $X_t$.

\begin{theorem}\label{thm1}
Let  
\begin{enumerate}

\item[(i)]
${\cal L}(X_0) \in {\cal P}_2$.

\item[(ii)] the functions $b$ and $\sigma$ be Borel bounded and globally Lipschitz in $y$???,  
\item[(iii)]  
$\Sigma\Sigma^*[t,x,\mu]$ be continuous in $(t,x)$ for any $\mu$ uniformly,
%quadratic, symmetric, and 
and also uniformly nondegenerate:
\begin{equation}\label{si}
\inf\limits_{s\ge 0, x\in \mathbb R^d,\mu\in {\cal P}_2} \inf\limits_{|\lambda|=1}
\lambda^*\Sigma[s,x,\mu]\Sigma[s,x,\mu]^*\lambda >0,
%\lambda^*\sigma(x,y)\lambda >0.
\end{equation}

\item[(iv)]  $\Sigma[t,x,\mu]$ be uniformly (globally) Lipschitz with respect to the  $x$ variable, and 

\item[(v)] 
the vector function $B[t,x,\mu]$ admit a  
Dini type condition with respect to  $x$
%, i.e., %there exists $\epsilon>0$ such that 
\begin{equation}\label{Dinib}
\int_0^1 \frac{\rho_B(r)}{r}\,dr < \infty. 
\end{equation}
\end{enumerate}
Then 
%{\color{blue}under the assumption of (\ref{m4})} {\color{red}(nado ili net?)} 
solution of the equation (\ref{e1}) is  pathwise unique and, hence, it is  strong.  
%which is pathwise and weakly unique.
\end{theorem}

%(d=1??)

\begin{remark}
The condition (ii) on {\em boundedness} of coefficients can be relaxed; we use it just for an easier reference about weak existence as well as for an easier reference on parabolic equations. The ultimate goal to estaiblish strong uniqueness  in dimension $d>1$ under just boundedness and Borel measurability of the drift $b$ still seems to remain an open problem for the case of the diffusion depending on the measure variable; if diffusion does not depend on the measure variable, see, for example, \cite[Theorem 3]{MV} (even under the linear growth condition in $x$ variable). The classical assumption of this theory (\ref{e200}) can also be relaxed: all is required is, actually Lipschitz condition of $\Sigma$ and $B$ in the measure variable. We do not pursue it here so as to simplify the presentation as much as it is possiible.
\end{remark}

%(odnomernyj sluchaj?)
%\begin{theorem}\label{thm2}
%If $d=1$, 
%\end{theorem}

\section{Proof of Theorem \ref{thm1}}
%\noindent
%{\bf Proof of Theorem \ref{thm1}}\,
The idea is to use Zvonkin's transformation as in  \cite{Zvonkin} to tackle the dependence of coefficients on the measure variable. 
Existence of weak solutions under the variety of conditions which include those in theorem \ref{thm1} was  established by many authors, see, e.g., \cite{Chiang, Funaki, HammersleySiskaSzpruch,  SimaMehri2, MV}.
%; in any case, the reference to \cite[Theorem 1]{MV} suffices? (ne sovsem).

Suppose there are two solutions $(X^i_t, \mu^i_t)$, $i=1,2$, on the same probability space with the same Wiener process $(W_t)$. 
Let $T>0$ (it will be chosen small enough, see later), 
and let 
\begin{equation}\label{L}
L^i(t,x) = L[t,x,\mu^i_t] = \frac12\, \sum_{j,k}A_{jk}[t,x,\mu^i_t]\frac{\partial^2}{\partial x^i \partial x^k} + \sum_{j}B^j[t,x,\mu^i_t]\frac{\partial}{\partial x^j}, \; i=1,2, 
\end{equation}
where $A[t,x,\mu] = \Sigma\Sigma^*[t,x,\mu]$, 
and let $u(s,x)=(u^1(t,x),\ldots, u^d(t,x))$ be a vector-function of Sobolev $\bigcap_{p>1}W_{p, loc}^{1,2}$ continuous solutions of the parabolic equations
\[
u^k_t(t,x) + L^1_t u^k(t,x) = 0, \quad u(T, x) = x^k, \quad 1\le k\le d.
\]
Note that $\Sigma[t',x,\mu_t] = \mathbb E \sigma(t',x,X_t)$ is continuous in $t$ as a solution of some parabolic PDE for each $t',x$. 
Solution of this system exists and is unique in the class of continuous functions in $\bigcap_{p>1}W_{p, loc}^{1,2}$ with a moderate growth of the function $u$ itself (for example, no faster than any polynomial), cf. \cite[Theorem 5.4]{Solo65}, \cite{Ver82, Zvonkin} (even under a bit more relaxed conditions: the continuity assumption of the matrix $A$ in $t$ at this stage could be dropped).
Denote 
\[
Y^i_t := u(t,X^i_t), \quad 0\le t\le T. 
\]
If $T>0$ is small enough, then it is known \cite{Zvonkin} that the gradient of the vector-function $u$ is close to the identity operator $I_{d\times d}$; in particular, it is uniformly bounded by the sup-norm, 
\begin{equation}\label{gradbdd}
\sup_{t,x}\|\nabla_x u(t,x)\|_B <\infty, 
\end{equation}
and there exists $C>0$ such that 
\[
C^{-1}|Y^1_t - Y^2_t| \le |X^1_t - X^2_t| \le C^{}|Y^1_t - Y^2_t|, \quad 0\le t\le T.
\]
Ito--Krylov's formula applied to $u(t,X^i_t)$ componentwise  for $0\le t\le T$ reads, 
\begin{align*}
dY^1_t = du(t,X^1_t) = (u_t + L^1_t u)(t,X^1_t)dt
+ \Sigma^*[t,X^1_t, \mu^1_t]\nabla_x u(t,X^1_t)dW_t
 \\\\
= \Sigma^*[t,X^1_t, \mu^1_t]\nabla_x u(t,X^1_t)dW_t, 
\end{align*}
since $u_t + L^1_t u = 0$ in the Sobolev sense, in particular, in $L_{d+1}$ which guarantees the applicability of Ito--Krylov's formula.
Also, similarly, 
\begin{align*}
dY^2_t = du(t,X^2_t) = \left(u_t + L^2_t u\right)(t,X^2_t)dt 
+ \Sigma^*[t,X^2_t, \mu^2_t]\nabla_x u(t,X^2_t)dW_t
 \\\\
=  \left(u_t + L^1_t u\right)(t,X^2_t)dt + \left(L^2_t - L^1_t\right)u(t,X^2_t)dt  
+ \Sigma^*[t,X^2_t, \mu^2_t]\nabla_x u(t,X^2_t)dW_t
 \\\\
=  \left(L^2_t - L^1_t\right)u(t,X^2_t)dt  
+ \Sigma^*(t,X^2_t, \mu^2_t)\nabla_x u(t,X^2_t)dW_t.
\end{align*}
Also notice that 
\[
Y^1_0 = u(0,x_0) = Y^2_0. 
\]
So, the difference $Y^1_t - Y^2_t$ has a (vector-valued) stochastic differential, 
\begin{align*}
d(Y^1_t - Y^2_t) \\\\
= \left(\Sigma^*[t,X^1_t, \mu^1_t)\nabla_x u(t,X^1_t) 
- \Sigma^*[t,X^2_t, \mu^2_t]\nabla_x u(t,X^2_t)\right)dW_t
-  \left(L^2_t - L^1_t\right)u(t,X^2_t)dt
 \\\\
= \left(\Sigma^*[t,X^1_t, \mu^1_t]\nabla_x u(t,X^1_t) 
- \Sigma^*[t,X^2_t, \mu^2_t]\nabla_x u(t,X^2_t)\right)dW_t
 \\\\
- \frac12 \mbox{Tr} \left(\Sigma\Sigma^*[t,X^2_t,\mu^2_t] - \Sigma\Sigma^*[t,X^2_t,\mu^1_t]\right) u_{xx}(t,X^2_t)dt
 \\\\
- \left(B[t,X^2_t,\mu^2_t] - B[t,X^2_t,\mu^1_t]\right) u_{x}(t,X^2_t)dt.    
\end{align*}
Hence, 
\begin{align*}
(Y^1_t - Y^2_t)^2 
= 2 \int_0^t (Y^1_s - Y^2_s)\left(\Sigma^*[s,X^1_s, \mu^1_s]\nabla_x u(s,X^1_s) \right.
 \\
\left. - \Sigma^*[s,X^2_s, \mu^2_s]\nabla_x u(s,X^2_s)\right)dW_s
\\
\!-  \!\int_0^t \!(Y^1_s \!-\! Y^2_s)
u_{xx}(s,X^2_s) \mbox{Tr}\left(\Sigma\Sigma^*[s,X^2_s,\mu^2_s] - \Sigma\Sigma^*[s,X^2_s,\mu^1_s]\right) ds
 \\
-  2\int_0^t  (Y^1_s - Y^2_s)
u_{x}(s,X^2_s)\left(B[s,X^2_s,\mu^2_s] - B[s,X^2_s,\mu^1_s]\right) ds
 \\
+ \int_0^t \left(\sigma^*(s,X^1_s, \mu^1_s)\nabla_x u(s,X^1_s) 
- \sigma^*(s,X^2_s, \mu^2_s)\nabla_x u(s,X^2_s)\right)^2ds.
\end{align*}
As usual, $(Y^1_t - Y^2_t)^2$ is understood as an inner product equal to $|Y^1_t - Y^2_t|^2$; the value 
$$
(Y^1_s - Y^2_s) u_{x}(s,X^2_s)
\left(B[s,X^2_s,\mu^2_s] - B[s,X^2_s,\mu^1_s)\right)
$$
is also interpreted as an inner product, or, equivalently, as pairing
$$
\langle Y^1_s - Y^2_s, u_{x}(s,X^2_s)
\left(B[s,X^2_s,\mu^2_s] - B[s,X^2_s,\mu^1_s]\right) \rangle, 
$$
although, the latter notation will not be used in the sequel, and all other parts of the latter calcuus are understood similary. 
Now,  after taking expectation of the left and right sides, the stochastic integral disappears. Further, due to the Lipschitz condition on both coefficients $b$ and $\sigma$ and their boundedness, we have the following inequality, 
\begin{align*}
\left|(L^2_s - L^1_s)u(s,X^2_s)\right| 
% \\\\
\le 
\frac12 \left|u_{xx}(s,X^2_s)\mbox{Tr}\left(\Sigma\Sigma^*[s,X^2_s,\mu^2_s] - \Sigma\Sigma^*[s,X^2_s,\mu^1_s]\right)\right| 
 \\\\
+  \left|u_{x}(s,X^2_s)\left(B[s,X^2_s,\mu^2_s] - B[s,X^2_s,\mu^1_s]\right)\right|
 \\\\
%\le C |\int b(y)(\mu^2_s(dy)-\mu^1_s(dy))|
%CV_1(\mu^1_s, \mu^2_s)   
\le C  \,  {\mathbb E} |X^1_s - X^2_s|.
\end{align*}
Here the first bound holds 
due to the boundedness of both $u_x$ (see (\ref{gradbdd})) and $u_{xx}$, the latter being established in \cite{Ivanovich66}, see also \cite[Proof of Theorem 4]{Zvonkin}. The last inequality is a consequence of the property that for any Lipsctiz continuous function $f$, 
\begin{align*}
|F[\mu^1]-F[\mu^2]| = |{\mathbb E}f(\xi^1) - {\mathbb E}f(\xi^2)|
= |{\mathbb E}(f(\xi^1) - f(\xi^2))|
 %\\\\
\le L_f {\mathbb E} |\xi^1 - \xi^2|, 
\end{align*} 
where $L_f$ is the Lipschitz norm of the function $f$.

Moreover, 
%\ifshowup
\begin{align*}
\left|(Y^1_s - Y^2_s)\left(\nabla_x u(s,X^1_s) \Sigma^*[s,X^1_s, \mu^1_s] 
- \nabla_x u(s,X^2_s) \Sigma^*[s,X^2_s, \mu^2_s]\right)\right|
 \\\\
\le C  |(Y^1_s - Y^2_s)||(X^1_s - X^2_s)| 
+ C |(Y^1_s - Y^2_s)| {\mathbb E}|(X^1_s - X^2_s)| 
 \\\\
\le C |(Y^1_s - Y^2_s)|^2 + C|(Y^1_s - Y^2_s)| \sqrt{{\mathbb E}|(Y^1_s - Y^2_s)|^2}, 
\end{align*}
by the Cauchy - Bouniakovskii -- Schwarz inequality. %\footnote{I wonder why Victor  Bouniakovskii is so unknown in Europe? He did his BSc and PhD at Sorbonne (under Cauchy supervision), his MSc in Germany, and extended Cauchy inequality for finite sums to integrals more than 3 decades later on,  and, yet, nearly 3 decades earlier than Schwarz. At least four countries must be  proud of him (and only two are). } 
We obtain, of course, with new constants $C$ in each line, 
\begin{align*}
{\mathbb E}_{} |Y^1_t - Y^2_t|^2 
\le C {\mathbb E}\int_0^t |(Y^1_s - Y^2_s)|^2ds 
 \\\\
+ C{\mathbb E}\int_0^t |(Y^1_s - Y^2_s)| ({\mathbb E}|(Y^1_s - Y^2_s)|^2)^{1/2} ds 
+  C{\mathbb E}\int_0^t |(Y^1_s - Y^2_s)|^2ds
 \\\\
\le C \int_0^t  {\mathbb E}|(Y^1_s - Y^2_s)|^2 ds.  
\end{align*}
Recall that the value ${\mathbb E} |(Y^1_s - Y^2_s)|^2$ is uniformly bounded for $0\le s\le T$, if $T>0$ is small enough (this was a requirement so as to guarantee the boundedness of $u_{xx}$). 
Hence, by virtue of Gronwall's inequality, 
\begin{align*}
{\mathbb E} |Y^1_t - Y^2_t|^2 = 0 \;\; \Longrightarrow \;\;  
{\mathbb E} |X^1_t - X^2_t|^2 = 0, \quad 0\le t\le T.  
\end{align*}
In other words, the solution is pathwise unique until $T$.
Repeating for $T\le t\le 2T$, $2T \le t\le 3T$, etc., we obtain pathwise uniqueness on the half line $t\ge 0$ by induction. The theorem is proved.

\section*{Acknowledgements}
This study has been funded by the Russian Academic Excellence Project '5-100' 
%(Theorem \ref{thm1})
%and by the RFBR grant \mbox{17-01-00633$\_$a}.
and 
%funded
by the Russian Science Foundation project no. 17-11-01098.
% (Theorem \ref{thm2}).


\begin{thebibliography}{99}

\bibitem{Bahlali19}
%Khaled Bahlali, Mohamed Amine Mezerdi, Brahim Mezerdi,  
K. Bahlali, M.A. Mezerdi, B. Mezerdi, Stability of McKean-Vlasov stochastic differential equations and applications,  
%{\em arXiv:1902.03478 [math.PR]}
{\em Stochastics and Dynamics,} 20(01), 2050007,  2020,
https://www.worldscientific.com/doi/10.1142/S0219493720500070
%no pages in the journal reference


\bibitem{Paul-Eric}
P.-E. Chaudru de Raynal, 
Strong well-posedness of McKean-Vlasov stochastic differential equation with H\"older drift, 
%{\em https://arxiv.org/abs/1512.08096}
Stochastic Processes Appl., 
130(1), 79--107, 2020.



\bibitem{Chiang}
T.S. Chiang,
 McKean-Vlasov equations with discontinuous coefficients.
 {\em Soochow J. Math.}, 20(4), 507--526, 1994.
% Dedicated to the memory of Professor Tsing-Houa Teng.

\bibitem{Funaki}
T.~Funaki.
 A certain class of diffusion processes associated with nonlinear parabolic equations.
 {\em Z. Wahrsch. verw. Gebiete}, 67(3): 331--348, 1984.


\bibitem{HammersleySiskaSzpruch}
W. Hammersley, D. \v{S}i\v{s}ka, L.  Szpruch, 
McKean-Vlasov SDEs under measure dependent Lyapunov conditions, 
{\em https://arxiv.org/abs/1802.03974}


\bibitem{Ivanovich66}
%М. Д. И в а н о в и ч , О характере непрерывности решений параболических урав­нений второго порядка, Вестник МГУ, серия матем., No 4 (1966), 31--41.
M.D. Ivanovi\v{c}, On the nature of continuity of solutions of linear parabolic equations of the second order, 
{\em Vestnik Moskov. Univ. Ser. I Mat. Meh.,} 21(4),  1966, 31--41 (Russian, with English summary).
%Moscow University Math.  Bulletin (Vestnik Moskovskogo Universiteta, Seriya 1, Matematika, Mekhanika),  no. 4, 31--41, 1966 (in Russian).
%{\color{red}(And where about Ivanovich refers further? to have a look!)}

%\bibitem{BJ97}
%B.~Jourdain.
% Diffusions with a nonlinear irregular drift coefficient and probabilistic interpretation of generalized {B}urgers' equations. {\em ESAIM Probab. Statist.}, 1:339--355 (electronic), 1995/97.

%\bibitem{BJSM98}
%B.~Jourdain and S.~M{\'e}l{\'e}ard. Propagation of chaos and fluctuations for a moderate model with smooth initial data.
%{\em Ann. Inst. H. Poincar\'e Probab. Statist.}, 34(6):727--766,1998.


\bibitem{Kry}
N.V. Krylov.
 {\em Controlled diffusion processes},  Springer-Verlag, Berlin et al., 2nd edition, 2009.  
 
%%N.~V. Krylov, Introduction to the Theory of Random Processes, AMS, Providence, R.I., 2002.

\bibitem{McKean}
H.P. McKean.
 A class of Markov processes associated with nonlinear parabolic
  equations.
 {\em Proc. Nat. Acad. Sci. U.S.A.}, 56:1907--1911, 1966.

\bibitem{Lad}
O.A. Ladyzhenskaya, V.A. Solonnikov, N.N. Uraltseva, {\em Linear and quasilinear equations of parabolic type,} Providence, R.I., AMS,  1968.


\bibitem{MehriStannat}
S. Mehri, W. Stannat, 
%Weak Solutions to Vlasov-McKean Equations under Lyapunov-type Conditions, https://arxiv.org/abs/1901.07778
Weak solutions to Vlasov–McKean equations under Lyapunov-type conditions, 
{\em Stochastics and Dynamics,} 19(06), 1950042: 1-23,  2019, 
%Free Access
https://doi.org/10.1142/S0219493719500424

\bibitem{SimaMehri2}
%S. Mehri, M. Scheutzow, W. Stannat, B.Z. Zangeneh, Propagation of Chaos for Stochastic Spatially Structured Neuronal Networks with Delay driven by Jump Diffusions, https://arxiv.org/abs/1805.01654

%Mehri, Sima, Scheutzow, Michael, Stannat, Wilhelm, and Zangeneh, Bian Z.
S. Mehri, M. Scheutzow, W. Stannat, B.Z. Zangeneh, 
Propagation of chaos for stochastic spatially structured neuronal networks with delay driven by jump diffusions,
{\em Annals of Applied Probability,} 30(1), 175--207, 2020.



\bibitem{MV}
Yu.S. Mishura, A.Yu. Veretennikov, 
Existence and uniqueness theorems for solutions of McKean--Vlasov stochastic equations, 
{\em https://arxiv.org/abs/1603.02212}

\bibitem{Shi}
A.N. Shiryaev, 
{\em Probability,}
Springer-Verlag, New York, 2nd edition, 1996.


\bibitem{Solo65}
%В. А. Солонников, “О краевых задачах для линейных параболических систем дифференциальных уравнений общего вида”, Тр. МИАН СССР, 83 (1965),  3–163      ; 
V.A. Solonnikov, On boundary value problems for linear parabolic systems of differential equations of general form, {\em Proc. Steklov Inst. Math.,} 83, 1--184, 1965.

\bibitem{Sz}
A.-S. Sznitman.
 Topics in propagation of chaos.
 In {\em \'{E}cole d'\'{E}t\'e de Probabilit\'es de Saint-Flour
  {XIX}---1989}, {\em Lecture Notes in Math.}, vol. 1464,  165--251.   Springer, Berlin, 1991.


%\bibitem{Ver80}
%A.~Yu. {Veretennikov}.
%{On strong solutions and explicit formulas for solutions of stochastic integral equations.} 
%{\em {Math. USSR, Sb.}}, 39:387--403, 1981.


\bibitem{Ver82}
A.Yu. Veretennikov. Parabolic equations and It\^o's stochastic equations with coefficients discontinuous in the time variable. {\em Math. notes,} 1982, 31(4), 278-283.

\bibitem{Ver06}
A.Yu. Veretennikov.
 On ergodic measures for McKean-Vlasov stochastic equations.
 In {\em Monte Carlo and quasi-Monte Carlo methods 2004},  Niederreiter, H.; Talay, D. (Eds.), 
  471--486. Springer, Berlin, 2006.


\bibitem{YamadaWatanabe}
T.~Yamada and S.~Watanabe.
 On the uniqueness of solutions of stochastic differential equations.
 {\em J. Math. Kyoto Univ.}, 11, 155--167, 1971.

%\bibitem{KZ}
%A.~K. Zvonkin and N.~V. Krylov.
% Strong solutions of stochastic differential equations.
% In {\em Proceedings of the School and Seminar on the Theory of
%  Random Processes (Druskininkai, 1974), Part II (Russian)}, pages
%  9--88. Inst. Fiz. i Mat. Akad. Nauk Litovsk. SSR, Vilnius, 1975.

\bibitem{Zvonkin}
A.K. Zvonkin, A transformation of the phase space of a diffusion process that removes the drift, 
%Mat. Sb. (N.S.), 93(135):1 (1974),  129–149
{\em Math. USSR-Sb.,} 22(1), 129--149, 1974.





\end{thebibliography}
\end{document}